\documentclass{article}
\usepackage[dvips]{graphics}
\usepackage[dvips]{color}
\usepackage{epsfig}
\usepackage{amsthm}
\usepackage{amsmath}
\usepackage{amsfonts}
\usepackage{url}
\usepackage{algorithm}
\usepackage{algorithmic}
\theoremstyle{definition}
\newtheorem{theorem}{Theorem}
\newtheorem{lemma}{Lemma}
\newtheorem{corollary}{Corollary}

\newcommand{\R}{\mathbb{R}}

\newcommand{\G}{\Gamma}

\def\ds{\delta(S)}

\newcommand{\beq}[1]{\begin{equation} \label{eq:#1} }
\newcommand{\refeq}[1]{(\ref{eq:#1})}
\newcommand{\eeq}{\end{equation}}

\newcommand{\Aut}[1]{\mathrm{Aut}({#1})}

\newcommand{\Sy}[1]{Sym({#1})}
\newcommand{\M}[1]{\mathrm{Met}_{#1}}
\newcommand{\m}[1]{\mathrm{met}_{#1}}
\newcommand{\Cu}[1]{\mathrm{Cut}_{#1}}
\newcommand{\Hy}[1]{\mathrm{Hyp}_{#1}}
\newcommand{\cu}[1]{\mathrm{cut}_{#1}}
\newcommand{\cG}{\overline{G}}

\title{The isometries of the cut, metric\\ and hypermetric cones\footnote{This research was
partially supported by the Research School SOM of the University of Groningen} 
}
\urldef{\dimaphome}\url{http://www.thi.informatik.uni-frankfurt.de/~dima/}
\urldef{\dimapmail}\url{dima@thi.informatik.uni-frankfurt.de}
\urldef{\bghome}\url{http://www.eco.rug.nl/medewerk/goldengorin/}
\urldef{\bgmail}\url{b.goldengorin@eco.rug.nl}
\urldef{\adhome}\url{http://www.is.titech.ac.jp/~deza/}
\urldef{\admail}\url{deza@is.titech.ac.jp}

\author{Antoine Deza\thanks{
 Department of Mathematical and Computing Sciences,
Tokyo Institute of Technology,
Meguro-ku Ookayama, Tokyo 152-8552, Japan,
Email: \admail . Web: \adhome}
\and
Boris Goldengorin\thanks{
Department of Econometrics
and Operations Research, Faculty of Economic Sciences
University of Groningen, P.O. Box 800, 9700 AV Groningen
The Netherlands,
Email: \bgmail . Web: \bghome } \and
Dmitrii V. Pasechnik\thanks{Theoretische Informatik,
FB~15, 
Robert-Mayer Str. 11-15, Postfach 11 19 32,
60054 Frankfurt am Main, Germany. 
Supported by the DFG Grant {\sf SCHN-503/2-1}.
Email: \dimapmail . Web: \dimaphome .} \footnote{
{\em Corresponding author.}  }
}
\date{June 2, 2003}
\begin{document}
\maketitle
\begin{abstract}
We show that the symmetry groups of the cut cone $\Cu{n}$ and the
metric cone $\M{n}$ both consist of the isometries induced by the
permutations on $\{1,\dots,n\}$; that is,
$Is(\Cu{n})=Is(\M{n})\simeq\Sy{n}$ for $n\geq 5$.  For $n=4$ we have
$Is(\Cu{4})=Is(\M{4})\simeq\Sy{3}\times\Sy{4}$.  This result can be
extended to cones containing the cuts as extreme rays and for which
the triangle inequalities are facet-inducing.  For instance,
$Is(\Hy{n})\simeq\Sy{n}$ for $n\geq 5$, where $\Hy{n}$ denotes the
hypermetric cone.
\end{abstract}

\section{Introduction and Notation}
The ${n\choose 2}$-dimensional {\it cut cone} $\Cu{n}$ is usually
introduced as the conic hull of the incidence vectors of all the cuts
of the complete graph on $n$ nodes.  More precisely, given a subset
$S$ of $V_n=\{1,\dots,n\}$, the {\it cut} determined by $S$ consists
of the pairs $(i,j)$ of elements of $V_n$ such that exactly one of
$i$, $j$ is in $S$. By $\ds$ we denote both the cut and its incidence
vector in $\R^{n\choose 2}$; that is, $\ds_{ij}=1$ if exactly one of
$i$, $j$ is in $S$ and $0$ otherwise for $1\leq i<j\leq n$. By abuse
of notation, we use the term cut for both the cut itself and its
incidence vector, so $\ds_{ij}$ are considered as coordinates of a
point in $\R^{n\choose 2}$. The cut cone $\Cu{n}$ is the conic hull of
all $2^{n-1}-1$ nonzero cuts, and the {\it cut polytope} $\cu{n}$ is
the convex hull of all $2^{n-1}$ cuts.  The cut cone and one of its
relaxation - the {\it metric cone} $\M{n}$ - can also be defined in
terms of finite metric spaces in the following way.  For all $3$-sets
$\{i,j,k\}\subset\{1,\dots,n\}$, we consider the following
inequalities.
\begin{eqnarray}
x_{ij}-x_{ik}-x_{jk} &\leq & 0 ,\label{eq:1}\\
x_{ij}+x_{ik}+x_{jk} &\leq & 2 .\label{eq:2}
\end{eqnarray}
(\ref{eq:1}) induce the $3{n \choose 3}$ facets which define the
metric cone $\M{n}$. Then, bounding the latter by the ${n \choose 3}$
facets induced by (\ref{eq:2}) we obtain the {\it metric polytope}
$\m{n}$.  The facets defined by (\ref{eq:1}) (resp. by (\ref{eq:2}))
can be seen as {\it triangle} (resp. {\it perimeter}) {\it
inequalities} for distance $x_{ij}$ on $\{1,\dots,n\}$.  While the cut
cone is the conic hull of all, up to a constant multiple,
$\{0,1\}$-valued extreme rays of the metric cone, the cut polytope
$\cu{n}$ is the convex hull of all $\{0,1\}$-valued vertices of the
metric polytope.  The link with finite metric spaces is the following:
there is a natural $1-1$ correspondence between the elements of the
metric cone and all the semi-metrics on $n$ points, and the elements
of the cut cone correspond precisely to the semi-metrics on $n$ points
that are isometrically embeddable into some $l^m_1$,
see~\cite{ad82}. It is easy to check that such minimal $m$ is smaller
than or equal to ${n \choose 2}$.

One of the motivations for the study of these polyhedra comes from
their applications in combinatorial optimization, the most important
being the MAXCUT and multicommodity flow problems.
For instance, the {\em linear programming approach} to MAXCUT involves
considering cutting planes that are needed to be added to $\M{n}$ to
obtain $\Cu{n}$. These cutting planes define cones $C_n$ such that
$\Cu{n}\subseteq C_n\subseteq\M{n}$. Perhaps the most well-known example
of such a $C_n$ is the {\em hypermetric cone} $\Hy{n}$ which is defined
by facets induced by inequalities generalizing the triangle inequalities.
For a detailed study of those
polyhedra and their applications in combinatorial optimization we
refer to {\sc Deza} and {\sc Laurent}~\cite{DeLa97}.

\section{Main Result}
One important feature of the metric and cut polyhedra is their very
large symmetry group. We recall that the symmetry group $Is(P)$ of a
polyhedron $P$ is the group of isometries preserving $P$ and that an
isometry is a linear transformation preserving the Euclidean distance.
While the symmetry groups of the cut and metric polytopes are known,
the question whether the cut and metric cones admit no other isometry
than the ones induced by $\Sy{n}$ was open,
see~\cite{DeDu02,DGL91,Lau96}.  More precisely, for $n\geq 5$,
$Is(\m{n})=Is(\cu{n})$ and both are induced by permutations on
$V_n=\{1,\dots,n\}$ and {\it switching reflections by a cut} and, for
$n\geq 5$, we have $|Is(\m{n})|=2^{n-1}n!$, see~\cite{DGL91}. Given a
cut $\ds$, the switching reflection $r_{\ds}$ is defined by
$y=r_{\ds}(x)$ where $y_{ij}=1-x_{ij}$ if $(i,j)\in \ds$ and
$y_{ij}=x_{ij}$ otherwise.

The aim of this article is to show that 
$Is(\Cu{n})=Is(\M{n})\simeq\Sy{n}$ for $n\geq 5$ and 
$Is(\Cu{4})=Is(\M{4})\simeq\Sy{3}\times\Sy{4}$.
A part of Theorem~\ref{th:th} was conjectured in 
\cite{DeDu02} and was substantiated by
computer calculations of the automorphism group of the ridge graph
of $\M{n}$ for $n\leq 20$. We recall that {\em ridge graph} of
a polyhedra $C_n$ is the graph which vertices are the facets of
$C_n$, two facets being adjacent if and only if their intersection
is a face of codimension $2$ of $C_n$. In other words, the ridge
graph of $C_n$ is the skeleton of the dual of $C_n$.

\begin{theorem}\label{th:th}
The symmetry groups of the cones $\M{n}$ and $\Cu{n}$ are isomorphic 
to $\Sy{n}$ for $n\geq 5$ and to $\Sy{3}\times\Sy{4}$ for $n=4$.
\end{theorem}

\noindent
The proof of Theorem~\ref{th:th} is given in Section~\ref{proofs}.
In Section~\ref{metcone} we characterize $Is(\M{n})$, 
in Section~\ref{cutcone} we show that $Is(\Cu{n})=Is(\M{n})$
and in Section~\ref{midcone} we generalize Theorem~\ref{th:th}
in the following way.

\begin{theorem}\label{th:gen} Let $C_n$ be a cone satisfying
\begin{itemize}
\item[$(i)$] the cuts are extreme rays of $C_n$, 
\item[$(ii)$] the triangle inequalities are facet-inducing for $C_n$.
\end{itemize}
Then any isometry of $C_n$ is induced by a permutation on $\{1,\dots,n\}$.
\end{theorem}

\noindent
A cone $C_n$ satisfying the condition $(i)$ and
$(ii)$ of Theorem~\ref{th:gen} is cone satisfying
$\Cu{n}\subseteq C_n\subseteq\M{n}$. Apart from $\M{n}$ and $\Cu{n}$
themselves,
a well-known example of such a cone $C_n$ is the {\it hypermetric cone}
$\Hy{n}$ defined by the following {\it hypermetric inequalities}
(\ref{eq:3})
which generalize the triangle inequalities:
\begin{eqnarray}
\sum_{1\leq i<j\leq n}b_i b_j \: x_{ij} & \leq 0 & \mbox{ with } \sum_{i=1}^n b_i=1 \label{eq:3}
\end{eqnarray}

\noindent
We recall that $Is(\Hy{n})$ contains the isometries induced by the
permutations on $\{1,\dots,n\}$. 

\begin{corollary} 
The symmetry group of $\Hy{n}$ is isomorphic 
to $\Sy{n}$ for $n\geq 5$ and to $\Sy{3}\times\Sy{4}$ for $n=4$.
\end{corollary}

\section{Proofs}\label{proofs}
We first prove Theorem~\ref{th:th} for $\M{n}$ by showing that
its ridge graph $G_n$ for $n>4$ has the automorphism group $\Sy{n}$;
the symmetry group of $\M{4}$ is constructed directly.
We complete the proof of Theorem~\ref{th:th} by 
showing that $G_n$ is an induced
subgraph of the ridge graph of $\Cu{n}$ that is invariant under the
isometries of $\Cu{n}$. Finally, we prove Theorem~\ref{th:gen} by noticing
that $\Cu{n}$ can be replaced by any cone $C_n$ satisfying  
$\Cu{n}\subseteq C_n\subseteq\M{n}$. 
The group-theoretic notation used in the paper can be found e.g. 
in \cite{Cam99}.

\subsection{The group $Is(\M{n})$ for $n\geq 5$.}\label{metcone}	
\subsubsection{$Is(\M{n})$ for $n\geq 4$}
Note that the isometries act faithfully on the facets; that is, the only
isometry that stabilizes each facet of $\M{n}$ is the trivial one.
As each permutation on $V_n$ is an isometry,
in order to prove the statement for $n\geq 5$, it suffices to show that
the automorphism group $A=\Aut{G_n}$ of the ridge graph $G_n$ is isomorphic to
$\Sy{n}$. 

The facets of $\M{n}$ naturally correspond to $\{0,1,-1\}$ vectors of
length $\binom{n}{2}$ with one positive and two negative entries.
Two triangle facets $u$ and $v$ are adjacent in $G_n$, i.e.
intersect on a face of codimension 2,  
if and only if they are {\em non-conflicting}; 
that is, there is no position $ij$ such that
the corresponding entries $u_{ij}$ and $v_{ij}$ 
are nonzero and of opposite sign, see~\cite{DeDe94,DeLa97}.

As already observed in \cite{DeDe94}, instead of working with $G_n$,
it appears to be easier to work with its complement $\cG_n$
which has the same automorphism group.
Observe that if two vertices $u$ and $v$ are conflicting then they
have either exactly one nonzero entry in common, or all the
three nonzero entries in common. Indeed, if 
$u_{ik}u_{jk}u_{ij}v_{ij}\neq 0$ and $v_{ik'}v_{jk'}\neq 0$,
then either $k=k'$, and the latter holds, or $k\neq k'$, and the former
holds.

The subgraph
induced on the neighbours of a vertex $u$ of $\cG_n$ is isomorphic to
the disjoint union of $n-3$ hexagons on a common edge $(u',u'')$, see~\cite{DeDe94}.
It is easy to see that $u$, $u'$, and $u''$ all have the same
zero entries. From now on let us refer to such a type of 3-clique in 
$\cG_n$ as a {\em Triangle}.
The Triangles $\{u,u',u''\}$ form an orbit of $A$, 
as the number of common neighbours of an edge of Triangle is $n-2$, 
bigger than for the edges of another types, where it is just $2$.

\begin{figure}[htbp]
\begin{center}
 \input{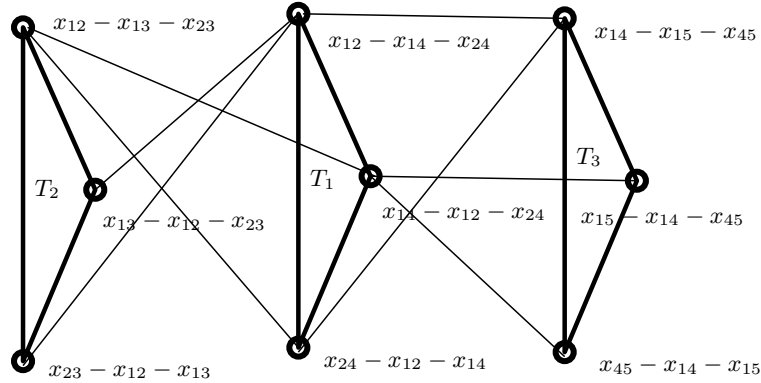}
\label{fig:tria}
\end{center}
\caption{Adjacencies between Triangles $T_1=\{1,2,4\}$, $T_2=\{1,2,3\}$,
and $T_3=\{1,4,5\}$.}
\end{figure}

Let us look at the edges between any two given Triangles $T_1$ and
$T_2$.  If they do not share a common nonzero entry then there are no
edges in between. Otherwise there are exactly 4 edges (see
Figure~\ref{fig:tria}), and, ignoring the edges of $T_1$ and $T_2$,
the subgraph induced on their 6 vertices is the disjoint union of two
2-paths.  This implies that every $g\in A$ stabilizing $T_1$ and $T_2$
either fixes $T_1$ pointwise, or induces an element of order $2$ on
the vertices of $T_1$. Let $T_3$ be a third Triangle having a common
entry with $T_1$, the different one than the common entry of $T_1$ and
$T_3$.  Then the 2-path between $T_1$ and $T_3$ with the middle point
in $T_1$ does not intersect the 2-path between $T_1$ and $T_2$ with
the middle point in $T_1$. Hence every $g\in A$ stabilizing $T_1$,
$T_2$ and $T_3$ fixes $T_1$ pointwise.  Therefore if $g\in A$
stabilizes all the Triangles, then $g$ is the identity; that is, $A$
acts faithfully on the set of Triangles.

Let us define the graph $\G_n$ on the Triangles, two Triangles are
adjacent if there is an edge of $\cG_n$ joining a vertex of the first
Triangle with a vertex of the second Triangle (so there are 4 edges
forming the disjoint union of 2-paths that join them).  To complete
the proof, it suffices to show that $\Aut{\G_n}\cong A$.
The latter is in fact not true for $n=6$, and we shall treat this case
separately.

Note that $\G_n$ is naturally 
isomorphic to the first subconstituent of the 
Johnson scheme $J(n,3)$, in other words, the graph with the vertex
set $\binom{V_n}{3}$, two vertices adjacent if the corresponding 3-subsets
intersect in a 2-subset. Automorphism groups of these graphs are described
e.g. in \cite{FKM94,KPR88}, and were known at least since
\cite{KlinPhD}. We give here a self-contained treatment, 
as the particular case we are dealing with is a simple one.
Note that $\G_5$ is the complement of the Petersen graph, and it is
well-known that its automorphism group is $\Sy{5}$.
This completes the proof for $n=5$.

For $n>5$, the graph $\G_n$ is distance-regular of diameter 3, 
see e.g. \cite{BCN} for this notion.
\begin{figure}[htbp]
\begin{center}
 \input{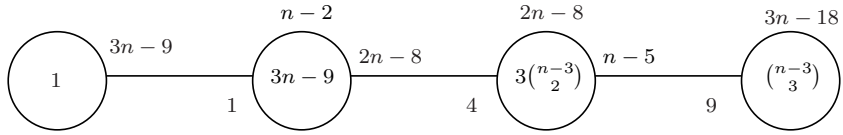}
\label{fig:jn3}
\end{center}
\caption{The distribution diagram of $\G_n$.}
\end{figure}

The subgraph $\Omega$ induced on the neighbourhood $\G_n(v)$ of a
vertex $v$ is isomorphic to the line graph of $K_{3,n-3}$.

If $n>6$ then the automorphism group of $\Omega$ is isomorphic to 
$\Sy{3}\times\Sy{n-3}$. It suffices to show that the
stabilizer of $v$ in $\Aut{\G_n}$ acts faithfully on $\G_n(v)$, as then it
will coincide with the one in $\Aut{\G_n}\cap\Sy{n}$.
Let $H$ be the kernel of the latter action.
As any vertex at distance 2 from $v$ is adjacent to exactly 4 vertices
in $\G_n(v)$, and any two different vertices at distance 2 from
$v$ are adjacent to different 4-sets of vertices in $\G_n(v)$, 
the vertices at distance 2 from $v$ are all fixed by $H$. 
In particular $\G_n(u)$ is fixed by $H$ for any $u\in\G_n(v)$.
Hence all the vertices at distance 2 from $u$ are fixed, too,
and $H=1$, as claimed. This completes the proof for $n>6$.

The graph $\G_6$ is a double antipodal cover of $K_{10}$, and
$\Aut{\G_6}\cong 2\times\Sy{6}$. We show that nevertheless only 
$\Sy{6}$ arises as $\Aut{\cG_6}$. Indeed, the normal subgroup $H=\Sy{2}$ 
interchanges simultaneously all the vertices of $\G_6$ at the maximal 
distance; they correspond to Triangles of $\cG_6$ with no common nonzero
coordinate. On the other hand, $H$ must act on $\cG_6$.
Observe that the pointwise stabilizer in $\Sy{6}$
of the facet $t_1=(x_{12},-x_{13},-x_{23})$, that belongs to a
Triangle
$T_1$, acts transitively on
the three facets with coordinates $x_{45}$, $x_{46}$ and $x_{56}$
forming a Triangle $T_2$.
On the other hand $H$ must map $t_1$ to one of these three latter
facets, as it must interchange $T_1$ and $T_2$.
This means that $H$ does not commute with the action of $\Sy{6}$ on 
the vertices of $\cG_6$, contradicting the assumption that $H\times\Sy{6}$
acts there. This completes the proof in the remaining case $n=6$.

\subsubsection{$Is(\M{4})$}\label{metcone4}
Here we shall construct the group of isometries of $\M{4}$ as 
the reflection group (see e.g. \cite{Hu90}) 
$A_2\times A_3\cong\Sy{3}\times\Sy{4}$. 
Let us recall the definition of a reflection $s(\alpha)$ with respect
to the hyperplane orthogonal to the vector $0\neq\alpha\in V$, 
with $V$ the Euclidean with the scalar product $\langle,\rangle$.
\beq{refl}
s(\alpha)v=v-2
\frac{\langle v,\alpha\rangle}{\langle \alpha,\alpha\rangle}\alpha,\quad
v\in V.
\eeq
Since $\M{4}=\Cu{4}$, the 7 extreme rays of 
$\M{4}$ are just the 7 nonzero cuts of the graph $K_4$.
In other words, a cut $V_1\cup V_2$ of $K_4$, i.e., a partition of the
vertices of $K_4$ into two parts $V_1$ and $V_2$, 
corresponds to the semimetric $d$ satisfying $d(x,y)=1$ for
$x$ and $y$ from different parts of the cut, and $d(x,y)=0$
otherwise. Writing the above-the-main-diagonal contents of $4\times 4$
symmetric matrices as vectors, the cuts $r_1$,\dots,$r_7$ are as follows:
\beq{rays}
\begin{matrix}
r_1=&(0& 0& 1& 0& 1& 1)\\
r_2=&(0& 1& 1& 1& 1& 0)\\
r_3=&(0& 1& 0& 1& 0& 1)\\
r_4=&(1& 0& 1& 1& 0& 1),\\
r_5=&(1& 1& 0& 0& 1& 1)\\
r_6=&(1& 0& 0& 1& 1& 0)\\
r_7=&(1& 1& 1& 0& 0& 0)\\
\end{matrix}\qquad
(r_4,r_5)=s(0, -1, 1, 1, -1, 0)
\eeq
Under the action of $\Sy{4}$, there are 2 orbits, one of 4 vectors
$r_1$, $r_3$, $r_6$, $r_7$ 
with 3 nonzero entries, and the other of the remaining 3 vectors.

We construct 5 reflections that generate the group 
$\Sy{3}\times\Sy{4}$
acting on  set \refeq{rays} of vectors. Each reflection is determined by the
hyperplane, generated by 5 linearly independent elements of \refeq{rays}.
More precisely, the reflection $s(\alpha)$ acting as the permutation
$(r_i,r_j)$ is given by an $\alpha\neq 0$ in the kernel
of the $5\times 6$ matrix 
$\left(\begin{matrix}r_{\ell_1}\\ \dots \\ r_{\ell_{5}}
\end{matrix}\right)$, where $\ell_j\neq i,j$.
Note that it is necessary that $r_i$ and $r_j$ lie in 
the same orbit of $\Sy{4}$.
For instance, $\alpha$ for $(r_4,r_5)$ is given in \refeq{rays}.

It is straightforward to check that indeed the 5 reflections just
described generate the group $A_2\times A_3$, and that they
act on the rays $r_i$'s in \refeq{rays}. Moreover, it suffices
to check that one of these reflections acts on the rays, as
together with the already present $\Sy{4}$ they generate the whole
group in question.

To complete the proof of Theorem~\ref{th:th} in this case, it suffices to refer to the
fact that the ridge graph $G_4$ of $\M{4}$ is isomorphic to the line
graph of $K_{3,4}$ (cf. \cite{DeDe94}), and thus the symmetry group cannot be
bigger than its automorphism group $\Sy{3}\times\Sy{4}$.

\subsection{The group $Is(\Cu{n})$ for $n\geq 4$}\label{cutcone}
As $\Cu{4}=\M{4}$, we can assume that $n>4$.
First, we remind that the maximal size facets of $\Cu{n}$
are the triangle facets given in~(\ref{eq:1}) and that a pair of
triangle facets are adjacent in the ridge graph of $\Cu{n}$ if and only if 
they are adjacent in the ridge graph of $\M{n}$.

\begin{lemma}[\cite{MR98c:90107}]\label{lem:trsizeC}
Any facet of $\Cu{n}$ contains at most  
$3\cdot 2^{n-3}-1$ extreme rays (cuts) with equality 
if and only if it is a triangle facet.
\end{lemma}

\begin{lemma}[\cite{DeDe94}]\label{lem:conflicC}
A pair of triangle facets of $\Cu{n}$ intersect on a face of 
codimension 2 if and only if they are non-conflicting.
\end{lemma}

\noindent
Lemma~\ref{lem:trsizeC} and Lemma~\ref{lem:conflicC}  
imply that the ridge graph of
$\M{n}$ is an induced subgraph in the ridge graph of $\Cu{n}$
that is invariant under any isometry of $\Cu{n}$.
Therefore we have $Is(\Cu{n})=Is(\M{n})$.

\subsection{The group 
$Is(C_n)$ for $\Cu{n}\subseteq C_n\subseteq\M{n}$}
\label{midcone}
Let $C_n$ be a cone having, among others, the cuts as extreme rays
and for which, among others, the triangle inequalities as facet-inducing.
As $\Cu{4}=\M{4}$, we can assume that $n>4$.
In the same way as for Lemma~\ref{lem:trsizeC}, Lemma~\ref{lem:trsizeCC}
can be directly deduced from a similar statement for the cut polytope $\cu{n}$,
see~\cite{MR96g:52022} and also \cite[Proposition~26.3.12]{DeLa97}. 

\begin{lemma}\label{lem:trsizeCC}
Let $C_n$ be a cone satisfying $\Cu{n}\subseteq C_n\subseteq\M{n}$,
any facet of $C_n$ contains at most $3\cdot 2^{n-3}-1$ cuts with equality 
if and only if it is a triangle facet.\qed
\end{lemma}

\noindent
Since $\Cu{n}\subseteq C_n\subseteq\M{n}$ and triangle facets are
adjacent in the ridge graphs of both $\Cu{n}$ and $\M{n}$
if and only if they are non-conflicting, we have the following.
\begin{lemma}\label{lem:conflicCC}
Let $C_n$ be a cone satisfying $\Cu{n}\subseteq C_n\subseteq\M{n}$,
a pair of triangle facets of $C_n$ intersect on a face of 
codimension 2 if and only if they are non-conflicting.\qed
\end{lemma}

\noindent
As in Section~\ref{cutcone}, 
Lemma~\ref{lem:trsizeCC} and Lemma~\ref{lem:conflicCC}  
imply that the ridge graph of
$\M{n}$ is an induced subgraph in the ridge graph of $C_n$
that is invariant under any isometry of $C_n$.
This completes the proof of Theorem~\ref{th:gen}.

\subsection*{Acknowledgment}
The authors thank 
Mikhail Klin and Monique Laurent 
for useful remarks and supplying relevant references.

{\small
\bibliography{../../tex/geom}
\bibliographystyle{abbrv}
}
\end{document}